\documentclass[12pt]{amsproc}

\RequirePackage{srcltx}

\usepackage{epstopdf}
\usepackage{booktabs}
\usepackage{amsmath}
\usepackage{amsthm}
\usepackage{amsfonts,amssymb}
\usepackage{color}
\usepackage{enumerate}
\usepackage{graphicx}

\newtheorem{thm}{Theorem}

\theoremstyle{remark}

\numberwithin{equation}{section}

\newcommand\NN{\mathbb{N}}
\newcommand\RR{\mathbb{R}}
\renewcommand\SS{\mathbb{S}}

\newcommand\proj{\operatorname{proj}}

\newcommand\esssup{\operatorname*{ess\,sup}}

\newcommand\Cdot{\mathop\cdot}

\newcommand\HH{\mathcal{H}}

\def\f{\frac}

\def\({\left(}
\def\){\right)}
\def\sa{\sigma}

\def\l{{\lambda}}

\def\o{{\omega}}
\def\s{{\sigma}}

\def\NN{{\mathbb N}}

\def\RR{{\mathbb R}}

\def\SS{{\mathbb S}}

\def\proj{\operatorname{proj}}

\def\s{\sa}

\def\HH{\mathcal{H}}

\begin{document}

\title{Nikolskii inequality for lacunary spherical polynomials}

\author{Feng Dai}
\address{F.~Dai, Department of Mathematical and Statistical Sciences\\
University of Alberta\\ Edmonton, Alberta T6G 2G1, Canada.}
\email{fdai@ualberta.ca}

\author{Dmitry~Gorbachev}
\address{D.~Gorbachev, Tula State University,
Department of Applied Mathematics and Computer Science, 300012 Tula, Russia}
\email{dvgmail@mail.ru}

\author{Sergey Tikhonov}
\address{S.~Tikhonov, Centre de Recerca Matem\`{a}tica\\
Campus de Bellaterra, Edifici~C 08193 Bellaterra (Barcelona), Spain; ICREA, Pg.
Llu\'{i}s Companys 23, 08010 Barcelona, Spain, and Universitat Aut\`{o}noma de
Barcelona.}
\email{stikhonov@crm.cat}

\thanks{F.~D. was supported
by NSERC Canada under the grant RGPIN 04702 Dai. D.~G. was supported by
the Russian Science Foundation under grant 18-11-00199.
S.~T. was partially supported by
MTM 2017-87409-P, 2017 SGR 358, and by the CERCA Programme of the Generalitat de Catalunya.}



\begin{abstract}
We prove that for $d\ge 2$, the asymptotic order of the usual
Nikolskii inequality on $\SS^d$ (also known as the reverse H\"{o}lder's inequality) can be significantly improved in many cases,
for lacunary spherical polynomials of the form $f=\sum_{j=0}^m f_{n_j}$ with
$f_{n_j}$ being a spherical harmonic of degree $n_j$ and $n_{j+1}-n_j\ge 3$. As
is well known, for $d=1$, the Nikolskii inequality for trigonometric
polynomials on the unit circle does not have such a phenomenon.
\end{abstract}

\date{\today}
\maketitle
\smallskip

\section{Introduction}

Let $\SS^{d}=\{x\in \mathbb{R}^{d+1}\colon |x|=1\}$ denote the unit sphere of
$\RR^{d+1}$ equipped with the usual surface Lebesgue measure $d\sigma(x)$, and $\omega_d$ the surface area 
of the sphere $\SS^{d}$; that is,
$\omega_{d}:=\s(\SS^d)=2\pi^{\frac{d+1}{2}}/\Gamma(\frac{d+1}{2})$. Here,
$|\Cdot|$ denotes the Euclidean norm of $\RR^{d+1}$. Given $0<p\le \infty$, we
denote by $L^p(\SS^{d})$ the usual Lebesgue $L^p$-space defined with respect to
the measure $d\sigma(x)$ on $\SS^{d}$, and
$\|\Cdot\|_p=\|\Cdot\|_{L^{p}(\SS^{d})}$ the quasi-norm of $L^p(\SS^{d})$; that
is,
\[
\|f\|_{p}:=\begin{cases}
\displaystyle\Bigl(\int_{\SS^{d}}|f(x)|^{p}\,d\sigma(x)\Bigr)^{1/p},&0<p<\infty,\\
\esssup_{x\in \SS^{d}}|f(x)|,&p=\infty.
\end{cases}
\]
 In what follows $c$, $C$ will denote 
 positive constants whose value may change with each occurrence.
 The notation $ A\asymp B$ means that 
$ c^{-1} A \le B
\le c A.$

Let $\Pi_n^d$ denote the space of all spherical polynomials of
degree at most $n$ on $\SS^{d}$ (i.e., restrictions on $\SS^{d}$ of polynomials
in $d+1$ variables of total degree at most $n$), and $\HH_n^d$ the space of all
spherical harmonics of degree $n$ on $\SS^{d}$. As is well known (see, e.g., \cite[Chap.~1]{BOOK}), both $\HH_n^d$ and $\Pi_n^d$ are finite
dimensional spaces with $\dim \HH_n^d\asymp n^{d-1}$ and $
\dim \Pi_n^d
\asymp n^{d}$.


The spaces $\HH_k^d$ are mutually orthogonal with respect to the inner product
of $L^2(\SS^{d})$, and the orthogonal projection $\proj_k$ of $L^2(\SS^{d})$
onto the space $\HH_k^d$ can be expressed as a spherical convolution:
\begin{equation}\label{1-1-0}
\proj_kf(x) =\frac{k+\lambda}{\lambda}\,\f1{\o_d}\int_{\SS^{d}} f(y) C_k^\lambda
(x\Cdot y)\,d\sigma(y),\quad x\in\SS^{d},\quad \l=\f{d-1}2,
\end{equation}
where the $C_k^\lambda$ denote the Gegenbauer polynomials as defined in \cite[Sec.~10.9]{BE53}.

The classical \textit{Nikolskii inequality} for spherical polynomials reads as
follows (see, e.g., \cite{Kam}):
\begin{equation}\label{Nikol}
\|f\|_q \le C_{d}n^{d\,(\frac 1p-\frac 1q)}\|f\|_p,\quad \forall\,f\in\Pi_n^d,\quad
0<p<q\le \infty.
\end{equation}

In \cite{pams}, continuing Sogge's investigations  \cite{sogge}, we obtained  the sharp asymptotic order
of the following Nikolskii  inequality for spherical harmonics $f_n$ of degree $n$, that is,
\begin{equation}\label{1-7--}
\|f_n\|_{q}\leq C n^{c(p,q)}\|f_n\|_{p},
 \quad
 \forall\,f_n\in\mathcal{H}_n^d,\quad
0<p<q\le \infty.
 \end{equation}
   In many cases, these  sharp estimates turn out to be  remarkably  better   than  the corresponding estimate for
   spherical polynomials (\ref{Nikol}). In particular,
   we have that
\begin{equation}\label{1-7}
\quad
 \|f_n\|_{q}\leq C n^{\f{d-1}2(\f 1p-\f1{q})} \|f_n\|_p,
\quad
 \forall\,f_n\in\mathcal{H}_n^d,\quad
  1\leq p\leq 2, \ \ \ p\leq q\leq p'.\end{equation}
Furthermore, this estimate is sharp.

These results, in particular, shows that there are no exponents  $0<p<q\le \infty$ such that the equivalence  $\|f_n\|_{q}\asymp \|f_n\|_p$ holds for any $f_n\in  \mathcal{H}_{n}^d$.  This in turn  implies that
no analogue of the following Zygmund theorem for  lacunary trigonometric series \cite{Zygmund}  for spherical polynomials:
{\it For any trigonometric series of the form
$$f(x)\sim\frac{a_0}2+\sum_{k=1}^\infty(a_k\cos n_kx+b_k\sin n_kx), \qquad \frac{n_{k+1}}{n_k}\ge \gamma>1,$$
one has
\begin{equation}\label{3-11}
 \|f\|_{L_p(\mathbb{T})}\asymp  \|f\|_{L_2(\mathbb{T})}, 
\qquad 0<p<\infty,
\end{equation}
where the  equivalent constants depend only on $p$ and $\gamma$.
}

In this paper, we prove that for $d\ge 2$, the asymptotic order of the Nikolskii inequality can be significantly improved when restricted
on a wide class of ``lacunary'' spherical polynomials, although the order is
sharp on the whole space of spherical polynomials.
To be precise, given
positive integers $\ell$, $m$ with $m\le n/\ell$, we denote by $\Pi_{n,m,\ell}^d$
the class of all spherical polynomials $f$ that can be represented in the form
\begin{equation}\label{3-2}
f=\sum_{j=0}^m f_{n_j},\quad f_{n_j}\in\mathcal{H}_{n_j}^d,\quad
j=0,1,\dots,m,
\end{equation}
for some sequence of nonnegative integers $\{n_j\}_{j=0}^m\subset [0,n]$ such
that $n_j-n_{j-1}\ge 2\ell+1$ for all $j=1,\dots, m$.
 Given $f\in \Pi_{n,m,\ell}^d$, we denote by $N_f$ the
largest integer $j$ for which $\proj_j f\neq 0$. Note that $\Pi_{n,m,\ell}^d$
is not a linear space.

Our main goal in this paper is to show

 \begin{thm}\label{lac}
   Let
 $(p,q)$ be a pair of exponents satisfying  either of the following  two conditions: (i) $0< p\le 1$ and $p\le q$ or (ii)  $1\le p\le 2$ and $p\le q\le p'$.   If  $f\in \Pi_{n,m,\ell}^d$, and  $d\ge 2$, then we have
\[
\|f\|_{q}\le C_d \bigl(n^{d-1-\ell_0}m\bigr)^{\frac 1p-\frac 1{q}}\|f\|_p\le C_d
n^{(d-\ell_0)(\frac 1p-\frac 1{q})}\|f\|_p,
\]
where $\ell_0=\min\{\ell, \frac{d-1}2\}$. In particular, if $\ell\ge \frac{d-1}2$,
then
\begin{equation}\label{3-2--}
\|f\|_{q}\le C_d \bigl(n^{\frac{d-1}2}m\bigr)^{\frac 1p-\frac 1{q}}\|f\|_p\le C_d
n^{\frac{d+1}2\,(\frac 1p-\frac 1{q})}\|f\|_p.
\end{equation}
\end{thm}


It is worth mentioning that  the asymptotic order of the Nikolskii exponent
for  ``lacunary'' spherical polynomials lies between the classical exponent provided by  (\ref{Nikol}) and the one for the spherical harmonics   given by (\ref{1-7--}).
In particular, we have  that for $\ell\ge \frac{d-1}2$,
 $$ \|f_n\|_{p'}\leq C n^{\f{d+1}2(\f 1p-\f1{p'})} \|f_n\|_p,\   \quad \forall f\in \Pi_{n,m,\ell}^d    \quad  1\leq p\leq 2.$$
It is also clear that inequality
(\ref{3-2--}) generalizes (\ref{1-7}).

Note that no improvement can be achieved in the order of the Nikolskii
inequality for similar ``lacunary'' trigonometric polynomials on the unit
circle, cf. (\ref{3-11}).

The   Nikolskii type   inequalities are closely related to the Remez type  inequalities in a very general setting, as was shown in  \cite[pp. 601-602]{teml}.  Moreover, these inequalities play a crucial role in establishing a  Sobolev-type embedding result for the Besov spaces: $B^r_{q}\big(L_{p}(\mathbb{S}^d)\big)\hookrightarrow L^{q}(\mathbb{S}^d)$,
(see \cite[Cor. 4]{hesse} and \cite[Sec.8]{dati}).
As a result, Theorem 1  can be applied to improve  the Remez type  inequalities
 as well as the limiting smoothness parameter $r=\frac{d+1}2\bigl(\frac1p-\frac1q\bigr)_+$ in  place of $r=d\bigl(\frac1p-\frac1q\bigr)_+$  for   ``lacunary'' spherical polynomials,  or
  ``lacunary'' spherical functions $f\in \bigcup_{n}\Pi_{n,m,\ell}^d$ with $\ell\ge \frac{d-1}2$.


\section{Proof of Theorem \ref{lac}
}\label{sec-lac-const}

We start with some useful definitions. Given $h\in\NN$, and a sequence
$\{a_n\}_{n=0}^\infty$ of real numbers, define (see, for instance, \cite{pams})
\[
\triangle_h a_n =a_n-a_{n+h},\quad \triangle_h^{\ell+1}=\triangle_h
\triangle_h^\ell,\quad \ell=1,2,\dots.
\]
Next, let
\[
R_n (\cos\theta):=\frac{C_n^{\lambda}(\cos\theta) }{C_n^{\lambda}(1)},\quad
\theta\in [0,\pi],
\]
denote the normalized Gegenbauer polynomial, and for a step $h\in \NN$, define
\[
\triangle_h^\ell R_n (\cos\theta):=\triangle_h^\ell a_n=\sum_{j=0}^{\ell}
(-1)^j \binom{\ell} j R_{n+hj} (\cos \theta),\quad \ell=1,2,\dots,\quad
n=0,1,\dots,
\]
with $a_n:= R_n(\cos\theta)$. Here and throughout, the difference operator in
$\triangle_h^\ell R_n (\cos\theta)$ is always acting on the integer $n$. In the
case when the step $h=1$, we have the following estimate
(\cite[Lemma~B.5.1]{BOOK}):
\begin{equation}\label{0-1}
\bigl|\triangle_1^{\ell} R_n (\cos\theta) \bigr|\le C \theta^\ell
(1+n\theta)^{-\frac{d-1}2},\quad \theta\in [0, \pi/2],\quad \ell\in \NN.
\end{equation}
On the other hand, however, the $\ell$-th order difference $\triangle_1^{\ell}
R_n (\cos\theta)$ with step $h=1$ does not provide a desirable upper estimate
when $\theta$ is close to $\pi$, and as will be seen in our later proof,
estimate \eqref{0-1} itself will not be enough for our purpose.

To overcome this difficulty, instead of the difference with step $1$, we
consider the $\ell$-th order difference $\triangle_2^{\ell} R_n (\cos\theta)$
with step $h=2$. Since $\triangle_2^\ell a_n =\sum_{j=0}^\ell\binom{\ell} j
\triangle_1^\ell a_{n+j},$ on one hand, \eqref{0-1} implies that
\[
\bigl|\triangle_2^{\ell} R_n (\cos\theta)
\bigr|\le C \theta^\ell (1+n\theta)^{-\frac{d-1}2},\quad \theta\in [0,\pi/2].
\]
On the other hand,
however, since
\[
\triangle_2^{\ell} R_n (\cos\theta) =\sum_{j=0}^{\ell} (-1)^j \binom{\ell} j
R_{n+2j} (\cos \theta),
\]
and since $R_{n+2j} (-z) = (-1)^n R_{n+2j} (z)$ (\cite[Sec.~10.9]{BE53}), we have $\triangle_2^{\ell}
R_n (\cos(\pi -\theta))=(-1)^n \triangle_2^{\ell} R_n (\cos\theta).$ It follows
that
\begin{equation}\label{0-2}
\bigl|\triangle_2^{\ell} R_n (\cos\theta) \bigr|\le C \begin{cases} \theta^\ell
(1+n\theta)^{-\frac{d-1}2}, & \theta\in [0, \pi/2],\\ (\pi-\theta)^{\ell} (1+n
(\pi-\theta))^{-\frac{d-1}2}, & \theta\in [\pi/2, \pi].\end{cases}
\end{equation}

By \eqref{1-1-0}, we obtain that for every $P\in \HH_n^d$,
\[
P(x)=c_n\,\f 1{\o_{d}}\int_{\SS^{d}} P(y) R_n (x\Cdot y) \, d\sigma(y),\quad x\in\SS^{d},
\]
where
\[
c_n:=\frac{\Gamma(\frac{d+1}2)}{2\pi^{(d+1)/2}}\,\frac{d+2n-1}{d+n-1}\,\frac{\Gamma(
d+n)}{\Gamma(n+1)\Gamma(d)}\asymp n^{d-1},
\]
and $x\Cdot y$ denotes the dot product of $x, y\in\RR^d$. Since $R_j
(x\,\cdot\,)\in\HH_j^d$ for any fixed $x\in\SS^{d}$, it follows by the
orthogonality of spherical harmonics that for any $P\in \mathcal{H}_n^d$, and
any $\ell\in\NN$,
\begin{align*}
 P(x)& = c_n \sum_{j=0}^\ell (-1)^j \binom {\ell}j\,\f 1{\o_{d}}\int_{\SS^{d}} P(y) R_{n+2j}
 (x\Cdot y)\, d\sigma(y) \\
 &=c_n\,\f 1{\o_{d}}\int_{\SS^{d}} P(y) \triangle_2^\ell R_n (x\Cdot y) \, d\sigma(y).
\end{align*}

By \eqref{0-2}, for any positive integer $\ell$,
\begin{equation}\label{2-4-1}
\bigl|\triangle_2^{\ell} R_n (\cos\theta) \bigr|\le C \min\{ n^{-\ell},
n^{-\frac{d-1}2}\}.
\end{equation}



Let $f\in\Pi_{n,m,\ell}^d$ be given in \eqref{3-2} with $n_m=n$. Define the
operator
\[
Tg(x):=\f 1{\o_{d}}\int_{\SS^{d}} g(y) H(x\Cdot y)\, d\sigma(y),\quad x\in\SS^{d},\quad
g\in L^1(\SS^{d}),
\]
where
\[
H(\cos \theta)=\sum_{k=0}^m c_{n_k} \sum_{j=0}^{\ell} (-1)^j \binom {\ell}j
R_{n_k+2j}(\cos \theta).
\]
Clearly,
\begin{equation*}
  Tg=\sum_{k=0}^m \sum_{j=0}^{\ell} (-1)^j \binom {\ell}j\frac{c_{n_k}}{c_{n_k+2j}}
\proj_{n_k+2j}(g),
\end{equation*}
and hence $Tf=f$.
Since, by \eqref{2-4-1},
\begin{equation*}
  \|H\|_\infty \le C\sum_{k=0}^m n_k^{d-1-\ell_0}\le C m n^{d-1-\ell_0},
\end{equation*}
it follows that
\[
\|T g\|_\infty \le C m n^{d-1-\ell_0}\|g\|_1,\quad \forall\,g\in L^1(\SS^{d}).
\]

On the other hand, by Plancherel's formula, we have
\[
\|Tg\|_2\le C \|g\|_2,\quad \forall\,g\in L^2(\SS^{d}).
\]
Thus, applying the Riesz-Thorin interpolation theorem, we deduce that for
$1\le p\le 2$,
\begin{equation*}
  \|T g\|_{p'}\le C (m n^{d-1-\ell_0})^{\frac 1p-\frac 1{p'}} \|g\|_p,\quad
  \forall\,g\in L^p(\SS^{d}).
\end{equation*}
Taking  $Tf=f$ we arrive at
\begin{equation}\label{vsp-}  \|f\|_{p'}\le C (m n^{d-1-\ell_0})^{\frac 1p-\frac 1{p'}} \|f\|_p,\quad
  1\le p\le 2.
\end{equation}
Further,  log-convexity of ${L^p}$ norms, namely
$\|f\|_q\le \|f\|_p^\theta\|f\|_{p'}^{1-\theta}$ with $\frac{\theta}p+\frac{1-\theta}{p'}=\frac1q$ and $0\le \theta\le 1$, implies
$$\|f\|_{q}\le C_d \bigl(m n^{d-1-\ell_0}\bigr)^{\frac 1p-\frac 1{q}}\|f\|_p,
$$
where $1\le p\le 2$ and $p\le q\le p'.$

To complete the proof we have to show that this inequality is valid for $0< p< 1$ and $p\le q$.
Let first $q=\infty$.
Using (\ref{vsp-}) with $p=1$, we have
\begin{eqnarray*}
\|f\|_{1}=
\||f|^{1-p}|f|^{p}\|_{1}&\le& \|f\|_\infty^{1-p} \| f\|_p^{p}
\\
&\le& C_d \bigl(m n^{d-1-\ell_0}\bigr)\|f\|_1  \|f\|_\infty^{-p} \|f\|_p^p.
\end{eqnarray*}
This yields that
\begin{equation}\label{vsp+} \|f\|_{\infty}\le C_d \bigl(m n^{d-1-\ell_0}\bigr)^{\frac 1p}\|f\|_p.
\end{equation}
If $q<\infty$, we write
$$
\|f\|_{q}=
\||f|^{1-\frac{p}q}|f|^\frac{p}q\|_{q}\le \|f\|_\infty^{1-\frac{p}q} \| f\|_p^{\frac{p}q}.
$$
Applying (\ref{vsp+}) implies
$$
\|f\|_{q}\le C_d \bigl(m n^{d-1-\ell_0}\bigr)^{\frac 1p-\frac 1q}\|f\|_p^{1-\frac pq}
 \| f\|_p^{\frac{p}q}=
C_d \bigl(m n^{d-1-\ell_0}\bigr)^{\frac 1p-\frac 1q}\|f\|_p,
$$
completing the proof.


\end{document}